\documentclass[11pt]{amsart}
\newtheorem{theorem}{Theorem}[section]
\newtheorem{proposition}{Proposition}[section]
\newtheorem{lemma}{Lemma}[section]
\newtheorem{corollary}{Corollary}[section]
\newtheorem{definition}{Definition}[section]

\theoremstyle{definition}
\newtheorem{remark}{Remark}[section]

\numberwithin{equation}{section}

\newcommand{\Nat}{{\mathbb N}}

\newcommand{\Real}{{\mathbb R}}
\newcommand{\Com}{{\mathbb C}}
\newcommand{\B}{{\mathcal B}}
\newcommand{\G}{{\mathcal G}}

\begin{document}
\title[General strong criterion for R.H.]{A general strong Nyman-Beurling criterion for the Riemann Hypothesis}
\author{Luis B\'{a}ez-Duarte\\ \ \\16 May 2005}
\date{16 May 2005, extracted from preprint dated 7 June 2004}

\begin{abstract}
For each $f:[0,\infty)\rightarrow\Com$ formally consider its co-Poisson or M\"{u}ntz transform $g(x)=\sum_{n\geq 1}f(nx)-\frac{1}{x}\int_0^\infty f(t)dt$. For certain $f$'s with both $f, g \in L_2(0,\infty)$ it is true that the Riemann hypothesis holds if and only if $f$ is in the $L_2$ closure of the vector space generated by the dilations $g(kx)$, $k\in\Nat$. Such is the case for example when $f=\chi_{(0,1]}$ where the above statement reduces to the strong Nyman criterion already established by the author. In this note we show that the necessity implication holds for any continuously differentiable function $f$ vanishing at infinity and satisfying $\int_0^\infty t|f'(t)|dt<\infty$. If in addition $f$ is of compact support then the sufficiency implication also holds true. It would be convenient to remove this compactness condition.
\end{abstract}

\maketitle

\section{Introduction}
\subsection{Preliminaries and notation}
The Riemann hypothesis shall be abbreviated as RH. We denote $L_p:=L_p(0,\infty)$, and, likewise we use $C_0$, $C_{00}$ to denote, respectively, the space of continuous functions on $[0,\infty)$ vanishing at infinity, and its subspace of compactly supported functions, whereas $C^1_0$, $C^1_{00}$ denote their corresponding subspaces of continuously differentiable functions.\\
\ \\ 
We let $\chi:=\chi_{(0,1]}$ be the characteristic function of $(0,1]$, and $\rho(x):=x-[x]$ the fractional part of $x$. We set $\rho_1(x):=\rho(1/x)$. Note that $\rho_1\in L_p$ for $1<p\leq\infty$. \\
\ \\
For each $\lambda>0$ the \emph{dilation} $K_\lambda$ defined on functions $f:[0,\infty)\rightarrow\Com$ is given by $K_\lambda f(x):=f(\lambda x)$.  $K_\lambda$ is a bounded operator on any $L_p$, $1\leq p \leq \infty$.\\
\ \\
For any $g:[0,\infty)\rightarrow\Com$ define $\B(g)$ to be the \emph{linear hull}  of $\{K_n g: n\in\Nat\}$. Obviously, if $g\in L_p$, then $\B(g)\subset L_p$.\\
\ \\
The \emph{Mellin transform} of $f$ at $s\in\Com$, denoted $\hat{f}(s)$, is defined by
$$
\hat{f}(s):=\int_0^\infty t^{s-1}f(t)dt.
$$
If $f\in L_p $ for all $p\in(1,\infty)$, then the above integral converges absolutely and uniformly on compacts subsets of the strip $0 < \Re s < 1$, so $\hat{f}(s)$ is an analytic function therein.\\
\ \\
The \emph{M\"{u}ntz operator} \footnote{M\"{u}ntz \cite{muntz} himself calls $x Pf(x)$ the \textit{Euler difference} (i.e., the difference between an integral in $(0,\infty)$ and its extended Riemann sums), whereas recently J. F. Burnol (\cite{burnol}, \cite{burnol2}) calls it the \textit{(M\"{u}ntz-)modified Poisson summation of $f$.}} $P$ is the linear operator defined formally on functions $f:[0,\infty)\rightarrow\Com$ by
\begin{equation}\label{muntzdef}
Pf(x):=\sum_{n=1}^\infty f(nx) - \frac{1}{x}\int_0^\infty f(t)dt.
\end{equation}
It is clear that $K_\lambda P = PK_\lambda$ for all $\lambda>0$, i.e., $P$ is an \emph{invariant} operator as defined in the author's paper \cite{baez0}. In section \ref{pproperties} we study some basic properties of $P$.

\subsection{A generalized strong Nyman-Beurling theorem}
We recall that the \textit{strong Nyman Beurling criterion} proved in \cite{baez1} is the following theorem which both simplified and strengthened the result in (\cite{nyman}, \cite{beurling}):

\begin{theorem}\label{strongnyman}
RH $\Longleftrightarrow$ $\chi\in\overline{\B(\rho_1)}^{L_2}$.
\end{theorem}

One obvious attempt to prove RH along these lines would be to tackle the Hilbert geometry problem of finding the distance from $\chi$ to the linear span of $\{K_a\rho_1: 1 \leq a \leq n\}$, and then to let $n\rightarrow\infty$. However, the scalar products $\int_0^\infty \rho(1/at)\rho(1/bt)dt$, given by Vasyunin's formulae (see \cite{vasyunin}, \cite{notes}) have  proved so far too complicated to compute the \emph{gramian determinants} involved. Nevertheless a generalization of Vasyunin's formula and a description in some depth of the \textit{autocorrelation function} $t \mapsto \int_0^\infty \rho(1/t)\rho(1/xt)dt$ is found in \cite{all}.   

Since it is verified trivially that
\begin{equation}\label{mchi}
P\chi=-\rho_1,
\end{equation}
Theorem \ref{strongnyman} can be rephrased as 
$$
\textnormal{RH} \Longleftrightarrow \chi\in\overline{\B(P\chi)}^{L_2}.
$$
We were thus led naturally to search for \emph{strong kernels} $f\in L_2$ having the property that $Pf\in L_2$, and such that the analogous equivalence implication holds true:
$$
\textnormal{RH} \Longleftrightarrow  f \in \overline{\B(Pf)}^{L_2}.
$$
It is quite simple to see from a slight adaptation of the proof of Theorem \ref{strongnyman} that \emph{step functions} with non-vanishing Mellin transform in the critical strip are strong kernels. We doubt at present, that this can be of any use. Furthermore, if $f$ is a strong kernel, then both $cf$ for all $c\not=0$, and $K_\lambda f$ for any $\lambda>0$ are also strong kernels. But we do \emph{not} know whether $f_1$, and $f_2$ being strong kernels implies that $f_1+f_2$ is a strong kernel, even if $\hat{f_1}(s)+\hat{f_2}(s)\not=0$ in the critical strip.  That is, the class of strong kernels is a vertexless cone in $L_2$, but it is not known whether it is a subspace (exclusion made of the point $0$). We rather turn our attention to a subspace of $L_2$ made up of differentiable kernels that might have interesting analytic properties, amongst them, but not exclusively, the possibility of simplifying the scalar products $\int_0^\infty Pf(ax)\overline{Pf(bx)}dx$.

\begin{definition}\label{good}
The class of \textbf{good kernels} $\G$ is the family of functions $f \in C^1_0 \cap L_1$ with $\int_0^\infty t|f'(t)|dt<\infty$.
\end{definition}

A simple important property of good kernels is the following:

\begin{lemma}\label{tft0}
If $f\in\G$ then $f(t)=o(1/t)$ at infinity, and $f\in L_p(0,\infty)$, 
($1 \leq p \leq \infty$).
\end{lemma}
\begin{proof}
For any finite $T>0$
$$
\int_0^T tf'(t)dt=Tf(T)-\int_0^T f(t)dt,
$$
which shows that there exists $\lim_{T\rightarrow\infty}Tf(T)=c$. If $c\not=0$ then we would obtain that $f\not\in L_1$. The $L_p$ statement is now obvious.
\end{proof}

\begin{remark}\label{remarkongood}
Not only are all good kernels $f\in L_p$ for $1\leq p\leq\infty$, it will transpire below from Proposition \ref{pfconvolution} that $Pf\in L_p$ for $1 < p\leq\infty$. Furthermore we shall see that $Pf\in C_0$. 
\end{remark}

Now we give a preliminary answer to our quest for non-trivial strong kernels. For any analytic function in the critical strip we shall denote
$$
Z(g):=\{s:\ g(s)=0,\ 1/2<\Re s<1\}.
$$ 

\begin{theorem}[General Strong Theorem]\label{generalstrongnyman}
The following implications are true.\\
Necessary condition for RH:
\begin{equation}\label{necessity}
\textnormal{RH}  \Longrightarrow 
\left(f\in\G \Longrightarrow f\in\overline{\B(Pf)}^{L_2}\right),
\end{equation}
Sufficiency condition for RH:
\begin{equation}\label{sufficiency}
\displaystyle \left(f\in\G\cap C_{00}\ \ \textnormal{\&}\ \ f\in\overline{\B(Pf)}^{L_2}\right)  
\Longrightarrow (Z(\zeta)\subset Z(\hat{f})).
\end{equation}
\end{theorem}
\noindent
The obvious corollary is
\begin{theorem}\label{generalstrongnyman2}
If $f\in\G\cap C_{00}$ and $Z(\hat{f})=\emptyset$, then
\begin{equation}\label{iff}
\textnormal{RH}  \Longleftrightarrow f\in\overline{\B(Pf)}^{L_2}.
\end{equation}
\end{theorem}

It is quite clear that we should like to \textbf{remove} the compact support condition in the sufficiency statement, but that has proved difficult to this moment.

In  Section \ref{pproperties} we shall first  establish some needed properties of the M\"{u}ntz operator, then we prove the general strong Theorem \ref{generalstrongnyman} in Section \ref{provemain}.

\section{Some properties of the M\"{u}ntz operator}\label{pproperties}
\subsection{Existence of $Pf$}
We first recall that the M\"{u}ntz operator is well defined for a space of functions containing $L_1$. This follows immediately from J. F. Burnol's work (Lemma 2.1 in \cite{burnol}). The neat proof of his lemma is worth repeating in detail.

\begin{lemma}[J.F. Burnol]\label{burnollemma}
If  $f\in L_1(\epsilon,\infty)$ for all $\epsilon>0$, then
$$
\int_0^\infty \sum_{n=1}^\infty |f(nx)|\frac{dx}{x}<\infty.
$$ 
\end{lemma}
An obvious consequence is
\begin{corollary}\label{poking}
$Pf(x)$ is well-defined a.e. for all $f\in\G$, and it is locally integrable.
\end{corollary}

\begin{proof}[Proof of Lemma \ref{burnollemma}]
Fix any $\epsilon>0$. Now use the monotone convergence twice in the following chain, first to take the sum out, then to put it back in:
\begin{eqnarray}\nonumber
\int_\epsilon^\infty \sum_{n=1}^\infty |f(nx)|\frac{dx}{x}
&=&
\sum_{n=1}^\infty \int_\epsilon^\infty |f(nx)| \frac{dx}{x}\\\nonumber
&=&
\sum_{n=1}^\infty \int_0^\infty \chi\left(\frac{\epsilon}{x}\right)|f(nx)| \frac{dx}{x}\\\nonumber
&=&
\sum_{n=1}^\infty \int_0^\infty \chi\left(\frac{n\epsilon}{x}\right)|f(x)| \frac{dx}{x}\\\nonumber
&=&
\int_0^\infty \sum_{n=1}^\infty\chi\left(\frac{n\epsilon}{x}\right)|f(x)| \frac{dx}{x}\\\nonumber
 &=&
\int_0^\infty \left[\frac{x}{\epsilon}\right]|f(x)| \frac{dx}{x}\\\nonumber
 &\leq&
 \frac{1}{\epsilon}\int_\epsilon^\infty|f(x)| dx. 
\end{eqnarray}
\end{proof}

\subsection{A convolution operator}

Let $f:[0,\infty)\rightarrow\Com$ be measurable. For every $\sigma\in\Real$ introduce the norm 
$$
\displaystyle \nu_\sigma(F) = \int_0^\infty t^\sigma |F(t)|dt,
$$
which is indeed the norm of the Banach space $L_1((0,\infty),t^\sigma dt)$.   We now define, formally at first, the operator $T_f$ by 

\begin{equation}\label{TF}
T_F G(x):=\int_0^\infty G(xt^{-1})F(t)dt.
\end{equation}
Note that $T_F$ is \emph{ invariant}, i.e., it commutes with all dilations $K_\lambda$. 

\begin{proposition}\label{Tbounded}
If $\nu_{\frac{1}{p}}(F) < \infty$ for some fixed $p\in[1,\infty]$, then the linear operator $T_F$  acts continuously from $L_p$ to $L_p$, with norm satisfying
\begin{equation}\label{Tnorm}
\|T_F\|_{L_p} \leq  \nu_{\frac{1}{p}}(F).
\end{equation}
\end{proposition}
\begin{remark}
A simple example of this operator is the following: for any $f\in\G$ we have
\begin{equation}\label{f_as_T}
f=-T_{f'}\chi.
\end{equation}
Here $\nu_1(f')<\infty$ and $\chi$ is in all $L_p$. This seemingly trivial identity plays a crucial r\^{o}le in the proof of the necessity criterion. 
\end{remark}

\begin{proof}
If $p=\infty$ the result is trivial. For $1\leq p <\infty$ let $G\in L_p$. Now we write the integral in (\ref{TF}) as a \textbf{true convolution} in the locally compact, multiplicative, abelian group $A=(0,\infty)^{\times}$ provided with Haar measure $t^{-1}dt$:

\begin{eqnarray}\nonumber
x^{\frac{1}{p}} T_F G (x)
&=&
\int_0^\infty (xt^{-1})^{\frac{1}{p}} G(xt^{-1})\hspace{1mm} t^{1+\frac{1}{p}}F(t) \frac{dt}{t}\\\nonumber
&=&
(\phi * \psi)(x),
\end{eqnarray}
where $\phi(t):=t^{\frac{1}{p}}G(t)$, $\psi(t):=t^{1+\frac{1}{p}}F(t)$, and $\phi\in L_p(A)$, $\psi\in L_1(A)$. Hence \emph{Young's inequality} $\|\phi * \psi\|_{L_p(A)} \leq \|\phi\|_{L_p(A)}\|\psi\|_{L_1(A)}$ becomes

\begin{equation}\label{normineq}
 \|T_F G\|_p \leq \nu_{\frac{1}{p}}(F) \|G\|_p. 
\end{equation}
\end{proof}

\subsection{$P$ as a convolution operator}
For some purposes the following representation doubles as an alternative definition of the M\"{u}ntz operator $P$.

\begin{proposition}\label{pfconvolution}
Let $f\in\G$, then
\begin{equation}\label{Mconvolution}
Pf(x)=T_{f'}\rho_1(x)=\int_0^\infty \rho_1(x t^{-1})f'(t)dt, 
\end{equation}
with 
\begin{equation}\label{normineqpf}
\|Pf\|_p \leq \nu_{\frac{1}{p}}(f')\|\rho_1\|_p, \ \ \ 1<p<\infty.
\end{equation}
A fortiori $Pf\in L_p$ for $1<p\leq\infty$.
\end{proposition}

\begin{remark}\label{pflp}
Several comments are in order about this representation for $Pf$: naturally the equation (\ref{Mconvolution}) is to be interpreted in the a.e. sense. It is easy to see furthermore that $P$ is not continuous relative to any $L_p(0,\infty)$-norm due to the troublesome presence of the derivative $f'$ in the convolution 
(\ref{Mconvolution}). On the other hand it is gratifying to see that for $f\in\G$ both $f, Pf \in L_p$ for $1<p\leq\infty$, very much so for $p=2$. The index $p=1$ fails because $\rho_1\not\in L_1$.

It is noteworthy too that one can prove that the representing integral for $Pf$ is a continuous function vanishing at infinity, whose value at the origin is $-\frac{1}{2}f(0)$, but, as that is of no use presently, we defer it to a future note.
\end{remark}

We can now prove the main proposition of this Section:

\begin{proof}[Proof of Proposition \ref{pfconvolution}]
 In view of (\ref{mchi}) we have for each finite $T>0$
\begin{eqnarray}\nonumber
\int_0^{T} \rho\left(\frac{t}{x}\right)f'(t)dt
&=&
\int_0^{T} \left(\frac{t}{x}-\sum_{n=1}^\infty\chi\left(\frac{nx}{t}\right)\right)f'(t)dt\\\label{limtotake}
&=&
\frac{1}{x}\int_0^{T} t f'(t) dt -
\int_0^{T} \sum_{n=1}^\infty\chi\left(\frac{nx}{t}\right)f'(t)dt,
\end{eqnarray}
where integrating by parts we get
\begin{equation}\label{lim1}
\int_0^{T} t f'(t) dt \rightarrow \int_0^\infty f(t) dt, \ \ \ T\rightarrow\infty
\end{equation}
since $Tf(T)\rightarrow 0$ (see Lemma \ref{tft0}). On the other hand the monotone convergence theorem yields

\begin{eqnarray}\nonumber
\int_0^{T}\sum_{n=1}^\infty \chi\left(\frac{nx}{t}\right)f'(t)dt
&=&
\sum_{n=1}^\infty \int_0^{T}\chi\left(\frac{nx}{t}\right)f'(t)dt\\\nonumber
&=&
\sum_{n\leq\frac{T}{x}} \int_{nx}^{T}f'(t)dt\\\nonumber
&=&
\left[\frac{T}{x}\right]f(T)-\sum_{n\leq\frac{T}{x}}f(nx).
\end{eqnarray}
Then using Burnol's Lemma \ref{burnollemma} we see that the above sum converges absolutely a.e. as $T\uparrow\infty$, and, since $T f(T)\rightarrow 0$, we get

\begin{equation}\label{lim2}
\int_0^{T}\sum_{n=1}^\infty \chi\left(\frac{nx}{t}\right)f'(t)dt
\rightarrow \sum_{n=1}^\infty f(nx), \ \ \ (T\rightarrow\infty).
\end{equation}
Therefore inserting (\ref{lim1}) and (\ref{lim2}) into (\ref{limtotake}) we obtain the desired representation (\ref{Mconvolution}) for 
$Pf(x)$.

Finally the norm inequality follows immediately from inequality (\ref{normineq}) of Proposition \ref{Tbounded}.
\end{proof}

\subsection{M\"{u}ntz's  formula}
It has been shown by M\"{u}ntz \cite{muntz} that for all $f\in C^1$ with 
$f(x)$, and $x f'(x)$ of order $x^{-\alpha}$ at infinity for some $\alpha>1$, one has
\begin{equation}\label{muntzformula}
\zeta(s)\hat{f}(s)=(Pf)^{\wedge}(s), \ \ \ (0<\Re s <1).
\end{equation}
It is interesting to remark that for the special \emph{non-smooth kernel} $f=\chi$, where we recall that $\rho_1=-P\chi$,  the formula is also valid, in fact it is the \emph{proto-M\"{u}ntz formula}, namely
\begin{equation}\label{rhomuntz}
\frac{\zeta(s)}{-s}=\int_0^\infty x^{s-1}\rho_1(x)dx, \ \ \ 0<\Re s <1, 
\end{equation}
which is the basis both, for M\"{u}ntz's original proof of (\ref{muntzformula}), as well as for the Nyman-Beurling theorem, and our earlier strengthened Theorem \ref{strongnyman}. M\"{u}ntz showed in \cite{muntz} that numerous, if not most proofs of the functional equation for the Riemann zeta-function are derived from (\ref{muntzformula}). We have given an even more general recipe that produces a potentially infinite number of such proofs in \cite{baez2}. 

We now present a somewhat more general version of  M\"{u}ntz formula for smooth $f$. Not \emph{the} more general version, of course, but certainly sufficient for our immediate purposes. There seems to be no absolutely most general formulation; however, J. F. Burnol has gone into this matter in great depth in \cite{burnol2}.

\begin{theorem}[M\"{u}ntz's formula]\label{mellinmuntz}
If $f\in\G$, then M\"{u}ntz's formula (\ref{muntzformula}) holds true.
\end{theorem}

\begin{proof}
Since $f$ and $Pf$ are in $L_p $ for $1<p<\infty$ (see Remark \ref{pflp}), then the Mellin transforms in M\"{u}ntz's formula (\ref{muntzformula}) are well-defined, absolutely convergent integrals in the strip $0 < \Re s < 1$. Do note also that $\int_0^\infty t^\sigma |f'(t)|dt < \infty$ for $0\leq \sigma \leq 1$. On the right-hand side of (\ref{muntzformula}) we substitute $Pf$ by its pseudo-convolution representation in (\ref{Mconvolution}) and employ Fubini's theorem together with (\ref{rhomuntz}) to obtain

\begin{eqnarray}\nonumber
\int_0^\infty x^{s-1}Pf(x)dx
&=&
\int_0^\infty x^{s-1}\int_0^\infty \rho\left(\frac{t}{x}\right)f'(t) dt dx\\\nonumber
&=&
\int_0^\infty \int_0^\infty x^{s-1}\rho\left(\frac{t}{x}\right)dx\  f'(t) dt \\\nonumber
&=&
\int_0^\infty \int_0^\infty (ut)^{s-1}\rho\left(\frac{1}{u}\right)du\  tf'(t) dt \\\nonumber
&=&
\int_0^\infty u^{s-1}\rho\left(\frac{1}{u}\right)du \int_0^\infty t^s f'(t) dt \\\nonumber
&=&
\frac{\zeta(s)}{-s}\int_0^\infty t^s f'(t) dt\\\nonumber
&=&
\zeta(s)\int_0^\infty t^{s-1}f(t)dt,
\end{eqnarray}
where the last integration by parts above uses the fact that $tf(t)\rightarrow 0$ as $t\rightarrow\infty$ (Lemma \ref{tft0}).

It remains to justify the interchange of integrations involved in the second equality of the above chain. To do this follow the same steps with the appropriate absolute values; indeed, let $\sigma=\Re s$, $0<\sigma<1$, then
\begin{eqnarray}\nonumber
\int_0^\infty \int_0^\infty x^{\sigma-1}\rho\left(\frac{t}{x}\right)|f'(t)| dt dx
&=&
\int_0^\infty \int_0^\infty x^{\sigma-1}\rho\left(\frac{t}{x}\right)dx\ |f'(t)| dt\\\nonumber
&=&
\int_0^\infty \int_0^\infty (ut)^{\sigma-1}\rho\left(\frac{1}{u}\right)du\ t|f'(t)| dt \\\nonumber
&=&
\frac{\zeta(\sigma)}{-\sigma} \int_0^\infty t^\sigma |f'(t)| dt<\infty.
\end{eqnarray}
This concludes the proof. \end{proof}

\section{The proof of the main Theorem \ref{generalstrongnyman}} \label{provemain}

\subsection{The necessity criterion}
No doubt it can qualify as daydreaming, but the necessity part of the theorem would seem like a good test to disprove RH if one holds on to the hope of finding an analytic kernel with simple scalar products. It is also worth remarking how simple the proof turns out to be in relation to the original Nyman-Beurling theorem, and even to that of the strengthened version. It is quite puzzling to note that no condition whatsoever is required of the Mellin transform of the kernel, which is quite surprising, as the Mellin transform $\hat{\chi}(s)=\frac{1}{s}$ plays a very important r\^{o}le in the proof of the necessity condition in our strong Nyman-Beurling Theorem \ref{strongnyman}. This bears reflecting on. 

\begin{proof}[Proof of the \textbf{necessity} criterion (\ref{necessity})]
Assume RH is true and let $f$ be a good kernel. By the strengthened Nyman-Beurling Theorem \ref{strongnyman} there is a sequence $\{h_n:n\in\Nat\}\subset\B(\rho_1)$ such that
\begin{equation}
\|h_n+\chi\|_2\rightarrow 0, \ \ \ (n\rightarrow\infty).
\end{equation}
There is no loss of generality in writing

$$
h_n=\sum_{a=1}^n c_{n,a}K_a\rho_1.
$$
Now we define $H_n:=T_{f'}h_n$. Note that the invariance of $T_{f'}$ and the important representation $Pf=T_{f'}\rho_1$ give us

$$
H_n=\sum_{a=1}^n c_{n,a}K_a Pf,
$$
so that $H_n\in\B(Pf)$. But $-f=T_{f'}\chi$ as noted before in equation (\ref{f_as_T}), hence
\begin{equation}
\|H_n-f\|_2=\|T_{f'}(h_n+\chi)\|_2\rightarrow 0, \ \ \ (n\rightarrow\infty),
\end{equation}
since $T_{f'}$ is a bounded operator from $L_2$ to itself as was shown in Proposition \ref{Tbounded}. This yields $f\in\overline{\B(Pf)}^{L_2}$.  
\end{proof}

\subsection{The sufficiency criterion}
We now prove the sufficiency criterion in Theorem \ref{generalstrongnyman}. One should certainly like to remove the condition that the good kernel be of compact support. This would allow the consideration of analytic kernels which may indeed prove to yield more tractable scalar products.

\begin{remark}\label{lp}
One obviously expects a generalization to a sufficiency criterion for the non-vanishing of $\zeta(s)$ in the half-plane $\Re s>1/p$ in terms of $L_p$ approximations. But we shall not bother with it at present, particularly since we have neither worked out the necessity criterion in the $L_p$-case, nor have done this yet for the strengthened Nyman-Beurling Theorem \ref{strongnyman}.   
\end{remark}

\begin{proof}[Proof of the \textbf{sufficiency} criterion (\ref{sufficiency})]
Let $f\in\G\cap C_{00}$ with $f\in\overline{\B(Pf)}^{L_2}$. Let $\zeta(s)=0$  for some fixed $s$ with $\sigma=\Re s \in (\frac{1}{2},1)$. Our aim is to prove that $\hat{f}(s)=0$ . By hypothesis there is a sequence $h_n\in\B(Pf)$ such that $\|h_n-f\|_2\rightarrow 0$. Now apply the invariant operator $V=I-2K_2$ as follows. Define $F:=Vf$, and $H_n:=Vh_n$.  By invariance of $V$, each $H_n\in\B(PF)$. Clearly $F\in\G\cap C_{00}$, and $\int_0^\infty F(x)dx=0$. These two facts together imply that there is an $A>0$ such $PF(x)=0$ for $x>A$, and thus every $H_n(x)=0$ also vanishes for $A>0$. Further, by $L_2$-continuity of $V$, $\|H_n-F\|_2\rightarrow0$. But M\"{u}ntz's formula (\ref{muntzformula}) implies $H_n^{\wedge}(s)=0$, and therefore
$$
-\hat{F}(s)=\int_0^A x^{s-1}(H_n(x)-F(x))dx.
$$
 But these integrals tend to zero by Schwarz's inequality since $x^{s-1}$ is square integrable in $(0,A)$. Hence $0=\hat{F}(s)=(1-2^{1-s})\hat{f}(s)$, so $\hat{f}(s)=0$.
\end{proof}

\bibliographystyle{amsplain}

\begin{thebibliography}{10}


\bibitem{baez0}
L. B\'{a}ez-Duarte, \textit{A class of invariant unitary operators} Adv. Math. \textbf{144 }(1999), no.1, 1-12.

\bibitem{notes}
L. B\'{a}ez-Duarte, M. Balazard, B. Landreau et E. Saias, \textit{Notes sur la fonction $\zeta$ de Riemann, 3}, Adv. Math. \textbf{149} (2000) 130-144.

\bibitem{baez2}
L. B\'{a}ez-Duarte, \textit{A general statement of the functional equation for the Riemann zeta-function} , Thoughts on the Riemann hypothesis 27, preprint, 2002.

\bibitem{baez1}
L. B\'{a}ez-Duarte, \textit{A strenghtening of the Nyman-Beurling criterion for the Riemann hypothesis} Rendiconti Accad. Lincei, \textbf{23} (2003) 5-11.

\bibitem{all}
L. B\'{a}ez-Duarte, M. Balazard, B. Landreau et E. Saias, \textit{\'{E}tude de l'autocorrelation multiplicative de la fonction ``partie fractionnaire"}, to appear in the Ramanujan Journal, (2004).

\bibitem{beurling}
A. Beurling, \textit{A closure problem related to the Riemann zeta-function}, Proc. Nat. Acad. Sci. U.S.A. \textbf{41} (1955), 312-314.

\bibitem{burnol}
J.-F. Burnol, \textit{On Fourier and Zeta(s)}, Habilitationsschrift, Universit\'{e} de Nice, (2001).

\bibitem{burnol2}
J.-F. Burnol, \textit{Co-Poisson intertwining: distribution and function theoretic aspects}, Universit\'{e} de Nice, (2002) pre-print.

\bibitem{muntz}
C. H. M\"{u}ntz, \textit{Beziehungen der Riemannschen $\zeta$-Funktion zu willk\"{u}rlichen reellen Funktionen}, Mat. Tidsskrift, B (1922), 39-47.

\bibitem{nyman}
B. Nyman, \textit{On the One-Dimensional Translation Group and Semi-Group in Certain Function Spaces}, Thesis, University of Uppsala, 1950, 55p.

\bibitem{vasyunin}
V. I. Vasyunin, \textit{On a biorthogonal system associated with the Riemann hypothesis} (in Russian), Algebra i Analiz \textbf{7} (1995), no. 3, 118-135, translation in St. Petersburg Math. J. \textbf{7} (1966), no. 3, 405-419.

\end{thebibliography}

\ \\
\ \\
\noindent Luis B\'{a}ez-Duarte\\
Departamento de Matem\'{a}ticas\\
Instituto Venezolano de Investigaciones Cient\'{\i}ficas\\
Apartado 21827, Caracas 1020-A\\
Venezuela\\
\email{lbaezd@cantv.net}

 \end{document}